\documentclass[a4paper,12pt]{article}

\usepackage{amssymb,epsf,graphics,amsmath}

\newtheorem{theo}{Theorem}[section]

\newtheorem{defi}{Definition}

\newtheorem{corollary}[theo]{Corollary}

\newenvironment{proof}{{\bf Proof }}{\hfill $\Box$}

\newenvironment{rmq}{{\bf Remark }}{}

\newcommand{\C}{\mathbb{C}}
\newcommand{\N}{\mathbb{N}}
\newcommand{\R}{\mathbb{R}}

\newcommand{\Bb}{\mathcal{B}}

\newcommand{\Hh}{\mathcal{H}}
\newcommand{\Kk}{\mathcal{K}}

\newcommand{\Tr}{\mathrm{Tr}}

\begin{document}


\title{\sc  A Lindblad Model for a Spin Chain\\Coupled to Heat Baths}

\author{Ameur DHAHRI \\\vspace{-2mm}\scriptsize Ceremade, UMR CNRS 7534\\\vspace{-2mm}\scriptsize Universit\'e Paris Dauphine\\\vspace{-2mm}\scriptsize place de Lattre de Tassigny\\\vspace{-2mm}\scriptsize 75775 Paris C\'edex 16\\\vspace{-2mm}\scriptsize France \\\vspace{-2mm}\scriptsize e-mail: dhahri@ceremade.dauphine.fr}

\date{}

\maketitle
\begin{abstract}
  We study a XY model which consists of a spin  chain coupled to heat baths. We give a repeated quantum interaction Hamiltonian describing this model. We compute the explicit form of the associated Lindblad generator in the case of the spin chain coupled to one, two and several heat baths. We further study the properties of quantum master equation such as approach to equilibrium, local equilibrium states, entropy production and quantum detailed balance condition. 
\end{abstract}
\section{Introduction}
The object of the quantum theory of open systems is to study the interaction of a quantum systems with very large ones. There are two different approaches which have usually been considered by physicists as well as mathematicians: The Hamiltonian and the Markovian approaches. 

The Hamiltonian approach consists of studying the reversible evolution of a small system in interaction with an exterior system and its main tools are: modular theory, $W^*$-dynamical system, Liouvillean...

The Markovian approach consists of studying the irreversible evolution of these systems in interaction picture. The interaction between the two systems is described by a quantum stochastic differential equation (quantum Langevin equation), a Lindblad generator (or Lindbladian) which is the generator of quantum Markovian semigroup... 

It is well-known that any quantum Markovian semigroup dilate a quantum stochastic differential equation in the sense of Hudson-Parthasarathy (cf [HP]). Moreover, its Lindblad generator is used for guessing the quantum master equation which allow to study the physical properties of a quantum system in interaction with a quantum field: quantum decoherence, approach to equilibrium, quantum detailed balance condition ...(cf [L], [Dav], [F],...). In the literature, in order to explicit the form of a Lindblad generator, we use the weak coupling limit which describes the passage from the Hamiltonian approach to the Markovian one.

Recently, in [AP] the authors consider the setup of a small system having repeated interactions, for a short duration $h$, with elements of a sequence of identical quantum systems. They prove that for a good choice of the repeated quantum interaction Hamiltonian, we get the explicit form for the associated Lindblad generator.

Here, we study a XY model which consists of a spin chain coupled to several heat baths. The heat bath is modeled by an infinite chain of identical spins. The full system is described by the means of a repeated quantum interaction Hamiltonian $H$ defined on the Hilbert space $\Hh_S\otimes \mathcal{\eta}$, where $\Hh_S=\otimes_{k=1}^N\,\C^2$ and $\mathcal{\eta}=\C^2$. After computing the Lindblad generator, we study the properties of the associated master equation when we discuss the case of the spin chain coupled to heat baths at the same inverse temperature $\beta$ and the case of distinct temperatures.

This paper is organized as follows. In Section 2 we compute the Lindblad generator describing the spin chain coupled to one and two heat baths at inverse temperatures $\beta$ and $\beta'$. In Section 3 we study the Markovian properties of the spin chain coupled to two heat baths. We give the explicit form of the stationary state $\rho^\beta$ of the associated master equation in the case of $\beta=\beta'$, this is proved in Subsection 3.1. The property of approach to equilibrium is studied in Section 3.2. The explicit form of the local equilibrium states is treated in Subsection 3.3. In Subsection 3.4 we compute the entropy production. If $\beta=\beta'$, we show that a quantum detailed balance condition is satisfied with respect to $\rho^\beta$, this is given in Section 3.5. Finally, in Section 4 we study the case of a spin chain coupled to $r$ $(2\leq r\leq N)$ heat baths.
\section{A Lindblad generator  for a spin  chain}
In this section, we give a repeated quantum interaction model associated to a spin chain coupled to one and two heat baths. We modelize the heat bath by an infinite chain of spins. We fuhrer give the GNS representation associated to the spin chain coupled to one piece (spin) of the heat bath (cf [AJ]). Finally, from [AP] we obtain the associated Lindblad generator. 
\subsection{Repeated quantum interaction model}
In this subsection, we present one of main results of repeated quantum interaction models. We refer the interested reader to [AP] for more details. 

Let us consider a small system $\Hh_0$ coupled with a piece of environment $\Hh$. The interaction between the two systems is described by an Hamiltonian $H$ which is defined on $\Hh_0\otimes \Hh$ and depending of time $h$. The associated unitary evolution during the interval $[0,h]$ of times is
$$\mathbb{L}=e^{-ihH}.$$
After the first interaction, we repeat this time coupling the same $\Hh_0$ with a new copy of $\Hh.$ Therefore, the sequence of the repeated  quantum interactions is described by the space 
$$\Hh_0\otimes\bigotimes_{\N^*}\Hh.$$
The unitary evolution of the small system in interaction picture with the $n$-th copy of $\Hh$, denoted by $\Hh_n$, is the operator $\mathbb{L}_n$ which acts as $\mathbb{L}$ on $\Hh_0\otimes\Hh_n$ and acts as the identity on copy of $\Hh$ different to $\Hh_n$. The discrete evolution equation describing this model is defined on $\Hh_0\otimes\bigotimes_{\N^*}\Hh$ as follows
\begin{equation}
\left\{
\begin{array}{lcc}
V_{n+1}=\mathbb{L}_{n+1}V_n \\
V_0=I
\end{array}
\right.
\end{equation}

Let $\{X_i\}_{\i\in\Lambda\cup\{0\}}$ be an orthonormal basis of $\Hh$ with $X_0=\Omega$ and let us consider the coefficients $(\mathbb{L}^i_j)_{i,j\in\Lambda\cup\{0\}}$ which are operators on $\Hh_0$ of the matrix representation of $\mathbb{L}$ in the basis $\{X_i\}_{i\in\Lambda\cup\{0\}}$. Then a natural basis of $\mathcal{B}(\Hh)$ is given by the family of operators $\{a^i_j,\;i,j\in\Lambda\cup\{0\}\}$ where
$$a_j^i(X_k)=\delta_{ik}X_j,\;\mbox{ for all },\,i,j,k\in\Lambda\cup\{0\}.$$
It is useful to notice that
$$\mathbb{L}=\sum_{i,j\in\Lambda\cup\{0\}}\,\mathbb{L}^i_j\otimes a_j^i.$$

Now let us put $\Psi=\otimes_{\N^*}\Omega$. Then from [AP], we have 
$$\langle\Psi,\;V_n^*(X\times I)V_n\rangle=L^n(X),\,\mbox{ for all }\,X\in\mathcal{B}(\Hh_0),$$
where $L(X)=\sum_{i\in\Lambda\cup\{0\}}\mathbb{L}^{0*}_iX\mathbb{L}^0_i$ is a completely positive map.

The following result is deduced from [AP].
\begin{theo}\label{a''}
Suppose that there exist operators $L^0_0,\,L^0_i,\;i\in\Lambda$ such that:
\begin{enumerate}
\item[i)] $\mathbb{L}_0^0=I+hL^0_0+o(h)$
\item[ii)] $\mathbb{L}_i^0=\sqrt{h}L_i^0+o(\sqrt{h}).$
\end{enumerate}
Then there exists a self-adjoint operator $H_0$ on $\Hh_0$ such that
$$\lim_{h\rightarrow0}\frac{L(X)-X}{h}=\mathcal{L}(X),\;\forall X\in\mathcal{B}(\Hh_0),$$
with
$$\mathcal{L}(X)=i[H_0,\,X]+\frac{1}{2}\sum_{i\in\Lambda}\,(2L^{0*}_iXL^0_i-XL_i^{0*}L_i^0-L_i^{0*}L_i^0X).$$
\end{theo}
\begin{proof}
Let $X\in\Bb(\Hh_0)$. Then we have
\begin{eqnarray}\label{L}
L(X)&=&\sum_{i\in\Lambda\cup\{0\}}\,U_i^{0*}XU_i^0\nonumber\\
   &=&X+h(L_0^{0*}X+XL_0^0+\sum_{i\in\Lambda}\,L_i^{0*}XL_i)+o(h).
\end{eqnarray}
Note that the operator $\mathbb{L}$ is unitary. This gives
$$\mathbb{L}_0^{0*}\mathbb{L}_0^0+\sum_{i\in\Lambda}\,\mathbb{L}_i^{0*}\mathbb{L}_i^0=I.$$
This implies that
$$I+h\bigl(L_0^{0*}+L_0^0+\sum_{i\in\Lambda}\,L_i^{0*}L_i^0\bigr)+o(h)=I.$$
Hence we obtain
$$L_0^{0*}+L_0^0=-\sum_{i\in\Lambda}\,L_i^{0*}L_i^0+o(1).$$
It follows that
$$L_0^0+\frac{1}{2}\,\sum_{i\in\Lambda}\,L_i^{0*}L_i^0=-\bigl(L_0^0+\frac{1}{2}\,\sum_{i\in\Lambda}\,L_i^{0*}L_i^0\bigr)^*+o(1).$$
Then there exists a self-adjoint operator $H_0$ on $\Hh_0$ such that
\begin{equation}\label{L0}
L_0^0+\frac{1}{2}\,\sum_{i\in\Lambda}\,L_i^{0*}L_i^0=-iH_0+o(1).
\end{equation}
Thus if we replace (\ref{L0}) in (\ref{L}), then the operator $L$ is written as
\begin{equation} \label{Lind}
L(X)=X+h\Bigl\{i[H_0,\;X]+\frac{1}{2}\sum_{i\in\Lambda}\,(2L_i^{0*}XL_i^0-XL_i^{0*}L_i^0-L_i^{0*}L_i^0X)\Bigr\}+o(h).
\end{equation}
This proves the above theorem.
\end{proof}
\subsection{Spin chains coupled to  one heat bath}
The system we consider here consists of $N$ spins, each of them described by the 2-dimensional Hilbert space $\mathcal{\eta}=\mathbb{C}^2$. Thus the Hilbert space of the spin chain is $\Hh_S=\otimes_{k=1}^N\,\C^2.$ The Hamiltonian of the system is given by
$$H_S=B\sum_{k=1}^N\,\sigma_z^{(k)}+\sum_{k=1}^{N-1}\,(J_x\sigma_x^{(k)}\otimes\sigma_x^{(k+1)}+J_y\sigma_y^{(k)}\otimes\sigma_y^{(k+1)}),$$
where
\begin{enumerate}
\item[]$\sigma_x=
\left(
\begin{array}{cc}
0 & 1\\
1 & 0
\end{array}
\right),\;\;
\sigma_y=\left(
\begin{array}{cc}
0 & -i\\
i & 0
\end{array}
\right),\;\;
\sigma_z=\left(
\begin{array}{cc}
1 & 0\\
0 & -1
\end{array}
\right)$
\end{enumerate}
and $B$ is a real number describing the influence of an external magnetic field in $z$-direction, while the interaction between nearest neighbors is described by $J_x,\,J_y\in\R$. 

In the following, we treat only the case that $J_x=J_y=1$ and we assume that the scalar $B$ is equal to 1. 

Now, let us consider a heat bath which is modeled by an infinite chain of 2-level atoms (spins). The quantum repeated interaction Hamiltonian of the system coupled at the first spin to this heat bath is written as
$$H=H_S\otimes I +I \otimes H_R+\frac{1}{\sqrt{h}}\,(\sigma_x^{(1)}\otimes\sigma_x+\sigma_y^{(1)}\otimes\sigma_y),$$
with $H_R=\sigma_z$.

Put
\begin{enumerate}
\item[]$\sigma_+=
\left(
\begin{array}{cc}
0 & 1\\
0 & 0
\end{array}
\right),\;\;
\sigma_-=\left(
\begin{array}{cc}
0 & 0\\
1 & 0
\end{array}
\right),\;\;
n_+=\left(
\begin{array}{cc}
1 & 0\\
0 & 0
\end{array}
\right),\;\;
n_-=\left(
\begin{array}{cc}
0 & 0\\
0 & 1
\end{array}
\right).$
\end{enumerate}

Let us consider the orthonormal basis $\{\Omega,\,X\}$ of $\C^2$, with
$$
\Omega=
\left(
\begin{array}{c}
1 \\
0 
\end{array}
\right),\;\;
X=\left(
\begin{array}{c}
0 \\
1 
\end{array}
\right).$$
Then in this basis we have
$$H=
\left(
\begin{array}{cc}
H_S+I & \frac{2}{\sqrt{h}}\,\sigma_-^{(1)}\\
\frac{2}{\sqrt{h}}\,\sigma_+^{(1)} & H_S-I
\end{array}
\right)$$
and the unitary evolution during the interval $[0,\,h]$ of times is given by
$$\mathbb{L}=
\left(
\begin{array}{cc}
I-ihI-ihH_S-2h\,\sigma_-^{(1)}\sigma_+^{(1)}+o(h^2) & -2i\sqrt{h}\,\sigma_-^{(1)}+o(h^{3/2})\\
-2i\sqrt{h}\,\sigma_+^{(1)}+o(h^{3/2}) & I+ihI-ihH_S-2h\,\sigma_+^{(1)}\sigma_-^{(1)}+o(h^2)
\end{array}
\right).$$

Let us define the scalar product on $M_2(\C)$ by
$$\langle A,\,B\rangle_\beta=\Tr(\rho_\beta\,A^*B),\;\,\forall A,\;B\in M_2(\C),$$
with 
$$\rho_\beta=\frac{e^{-\beta\sigma_z}}{\Tr(e^{-\beta\sigma_z})}=
\left(
\begin{array}{cc}
\beta_0 & 0\\
0 & \beta_1
\end{array}
\right)$$
is the equilibrium state at inverse temperature $\beta$ of a single spin.

Put
$$X_0=I,\;
X_1=\frac{1}{\sqrt{\beta_0}}\left(
\begin{array}{cc}
0 & 0\\
1 & 0
\end{array}
\right),\;
X_2=\frac{1}{\sqrt{\beta_1}}
\left(
\begin{array}{cc}
0 & 1\\
0 & 0
\end{array}
\right),\;
X_3=\frac{1}{\sqrt{\beta_0\beta_1}}
\left(
\begin{array}{cc}
\beta_1 & 0\\
0 & -\beta_0
\end{array}
\right).$$
It is clear that $\{X_0,X_1,X_2,X_3\}$ form an orthonormal basis of $M_2(\C)$ equipped by the scalar product $\langle\,,\,\rangle_\beta.$ 

The GNS representation of $(\C^2,\,\rho_\beta)$ is the triple $(\pi,\,\widetilde{\Hh},\,\Omega_R)$ where
\begin{enumerate}
\item[] $\bullet\;\Omega_R=I,$
\item[] $\bullet\;\widetilde{\Hh}=M_2(\C),$
\item[] $\bullet\;\pi:\,M_2(\C)\rightarrow\Bb(\widetilde{\Hh})$, such that $\pi(M)A=MA,\;\forall M,\;A\in M_2(\C).$
\end {enumerate}
Now we prove the following.
\begin{theo}\label{a'}
The Lindblad generator of the repeated quantum interaction model associated to the spin chain coupled to one heat bath at positive temperature $\beta^{-1}$ is given by
\begin{eqnarray*}
\mathcal{L}_1(X)=i\,[H_S,\,X]&+&2\beta_0\,[2\sigma_-^{(1)}X\sigma_+^{(1)}-\{n_-^{(1)},\,X\}]\\
                           &+&2\beta_1\,[2\sigma_+^{(1)}X\sigma_-^{(1)}-\{n_+^{(1)},\,X\}],
\end{eqnarray*}
for all $X\in \mathcal{B}(\Hh_S).$
\end{theo}
\begin{proof}
Set $\widetilde{\mathbb{L}}=\pi(\mathbb{L})$. In the basis $\{X_0,\,X_1,\,X_2,\,X_3\}$ we have
\begin{enumerate}
\item[]$\widetilde{\mathbb{L}}_0^0=I-ihH_S+ih(\beta_1-\beta_0)I-2h\,\beta_0 \,\sigma_-^{(1)}\sigma_+^{(1)}-2h\,\beta_1\,\sigma_+^{(1)}\sigma_-^{(1)}+o(h^2)$,
\item[]$\widetilde{\mathbb{L}}_1^0=-2i\sqrt{\beta_0}\,\sqrt{h}\,\sigma_+^{(1)}+o(h^{3/2}),$
\item[]$\widetilde{\mathbb{L}}_2^0=-2i\sqrt{\beta_1}\,\sqrt{h}\,\sigma_-^{(1)}+o(h^{3/2}),$
\item[]$\widetilde{\mathbb{L}}_3^0=o(h).$
\end{enumerate}
Set
\begin{enumerate}
\item[]$H_0=H_S+(\beta_0-\beta_1)I,$
\item[]$L_0^0=-iH_0-2\beta_0\,\sigma_+^{(1)}\sigma_-^{(1)}-2\beta_1\sigma_-^{(1)}\sigma_+^{(1)},$
\item[]$L_1^0=-2i\sqrt{\beta_0}\,\sigma_-^{(1)},$
\item[]$L_2^0=-2i\sqrt{\beta_1}\,\sigma_+^{(1)}.$
\end{enumerate}
Then it is clear that
$$L^0_0=-iH_0-\frac{1}{2}\sum_{i=1}^{2}L_i^0L_i^{0*}.$$
Therefore by using Theorem \ref{a''}, the result of the above theorem holds.\bigskip
\end{proof}

\begin{rmq}
Note that if N=1, then the Lindblad generator is written as
\begin{eqnarray*}
\mathcal{L}_1(X)=i\,[H_S,\,X]&+&2\beta_0\,[2\sigma_-X\sigma_+-\{n_-,\,X\}]\\
                           &+&2\beta_1\,[2\sigma_+X\sigma_--\{n_+,\,X\}],
\end{eqnarray*}
for all $X\in M_2(\C)$. This Lindbladian describes a two-levels atom in interaction with a heat bath. It is easy to show that the associated master equation has the properties of approach to equilibrium and the quantum detailed balance condition with respect to thermodynamical state of the spin at inverse temperature $\beta$ is satisfied. We refer the interested reader to [D] for more details. Moreover at zero temperature, that is $\beta=\infty$, we can prove in the same way at [D] that the associated quantum dynamical semigroup converges towards the equilibrium.
\end{rmq}
\subsection{Spin chains coupled to two heat baths}
In this subsection, we suppose that the spin chain is coupled to two heat baths respectively at the first and the $N$-th spin. Moreover, the two heat baths are supposed to be respectively at inverse temperatures $\beta$ and $\beta'$. Therefore the associated repeated quantum interaction Hamiltonian is of the form
$$H=H_S\otimes I +I \otimes H_R+\frac{1}{\sqrt{h}}\,(\sigma_x^{(1)}\otimes\sigma_x^{(L)}+\sigma_y^{(1)}\otimes\sigma_y^{(L)}+\sigma_x^{(N)}\otimes\sigma_x^{(R)}+\sigma_y^{(N)}\otimes\sigma_y^{(R)}),$$
where $(R)$, $(L)$ indicate the left and the right heat bath.

The proof of the following theorem is similar as the one of Theorem \ref{a'}.
\begin{theo}
The Lindblad generator associated to the spin chain coupled to two heat baths at inverse temperatures $\beta$ and $\beta'$ is given by
\begin{eqnarray*}
\mathcal{L}(X)=i\,[H_S,\,X]&+&2\beta_0\,[2\sigma_-^{(1)}X\sigma_+^{(1)}-\{n_-^{(1)},\,X\}]\\
                           &+&2\beta_1\,[2\sigma_+^{(1)}X\sigma_-^{(1)}-\{n_+^{(1)},\,X\}]\\
                           &+&2\beta_0'\,[2\sigma_-^{(N)}X\sigma_+^{(N)}-\{n_-^{(N)},\,X\}]\\
                           &+&2\beta_1'\,[2\sigma_+^{(N)}X\sigma_-^{(N)}-\{n_+^{(N)},\,X\}],
\end{eqnarray*}
for all $X\in \mathcal{B}(\Hh_S).$
\end{theo}
\section{Markovian properties of a spin chain coupled to two heat baths}

In this section, we describe the Markovian properties of the spin chain coupled to two heat baths at inverse temperatures $\beta$ and $\beta'$. We start by giving the associated quantum master equation. Moreover, we study the property of approach to equilibrium and we compute the local states. Finally for $\beta=\beta'$, we compute the entropy production and we study the quantum detailed balance condition.

Note that, in order to study the two last physical properties, we need to know explicitly the stationary state which is complicated to compute it in the case where $\beta\neq\beta'$.
\subsection{Quantum master equation}
For all density matrix $\rho\in\mathcal{B}(\Hh_S)$, the quantum master equation of the spin chain coupled to two heat baths at inverse temperatures $\beta$ and $\beta'$ is defined as
\begin{eqnarray*}
\mathcal{L}^*(\rho)=-i\,[H_S,\,\rho]&+&2\beta_0\,[2\sigma_+^{(1)}\rho\,\sigma_-^{(1)}-\{n_-^{(1)},\,\rho\}]\\
                           &+&2\beta_1\,[2\sigma_-^{(1)}\rho\,\sigma_+^{(1)}-\{n_+^{(1)},\,\rho\}]\\
                           &+&2\beta_0'\,[2\sigma_+^{(N)}\rho\,\sigma_-^{(N)}-\{n_-^{(N)},\,\rho\}]\\
                           &+&2\beta_1'\,[2\sigma_-^{(N)}\rho\,\sigma_+^{(N)}-\{n_+^{(N)},\,\rho\}].
\end{eqnarray*}
Note that dim $\Hh_S<\infty$. Therefore there exists a stationary state for the above master equation. Hence in order to prove the uniqueness of the equilibrium state we need the following theorem (cf [F]).
\begin{theo}\label{a}
Let $(\Theta_t)_t$ be a norm continuous quantum dynamical semigroup on $\Bb(\Kk)$ for some separable Hilbert space $\Kk$ whose generator $L$ is given by
\begin{eqnarray}\label{A}
L(A)=\sum_j\,V_j^*AV_j+KA+AK^*,
\end{eqnarray}
where $V_j\in\Bb(\Hh)$ and $K=iH-\frac{1}{2}\sum_j\,V_j^*V_j,\;H=H^*\in\Bb(\Kk)\;\,(L(I)=0)$.
Suppose that $(\Theta_t^*)_t$ has a stationary faithful state $\rho$. Then $\rho$ is the unique stationary state for $(\Theta_t^*)_t$ if and only if
$$\{H,\;V_j^*,\;V_j\}'=\C I.$$
\end{theo}

Put 
$$\rho^\beta=\rho_\beta\otimes\, ...\,\otimes\rho_\beta=\otimes_{i=1}^N\rho_\beta.$$
Now we prove the following.
\begin{theo}\label{same-temperature}
If $\beta=\beta'$, then $\rho^\beta$ is the unique faithful stationary state for the quantum dynamical semigroup $(e^{t\mathcal{L}^*})_{t\geq0}$.
\end{theo}
\begin{proof}
Note that it is straightforward to show that 
$$[\sigma_x\otimes\sigma_x+\sigma_y\otimes\sigma_y,\,\rho_\beta\otimes\rho_\beta]=0.$$ 
Therefore we get $[H_S,\,\rho^\beta]=0$. Moreover, if we note by $\mathcal{L}^*_d$ the dissipative part of $\mathcal{L}^*$, then it is easy to show that $\mathcal{L}^*_d(\rho^\beta)=0$. Hence we have $\mathcal{L}^*(\rho^\beta)=0$. Thus $\rho^\beta$ is a stationary state for the above master equation. 

Now let us consider an operator $A$ such that
$$A\in\{H_S,\,\sigma_+^{(1)},\,\sigma_-^{(1)},\,\sigma_-^{(N)},\,\sigma_+^{(N)}\}'.$$ 
In particular we have 
$$A\in\{\sigma_+^{(1)},\,\sigma_-^{(1)},\,\sigma_-^{(N)},\,\sigma_+^{(N)}\}'.$$
 This gives
 $$A=I^{(1)}\otimes A_1\otimes I^{(N)},$$
 where $A_1$ is an operator defined on $\otimes_{k=2}^{N-1}\,\C^2$. On the other hand, we have $A$ commutes with $H_S$. Therefore we get
\begin{eqnarray*}
&{}&\sigma_x^{(1)}\otimes[A_1,\,\sigma_x^{(2)}]\otimes I^{(N)}+\sigma_y^{(1)}\otimes[A_1,\,\sigma_y^{(2)}]\otimes I^{(N)}+\\
&{}&I^{(1)}\otimes[A_1,\,\sigma_x^{(N-1)}]\otimes\sigma_x^{(N)}+I^{(1)}\otimes[A_1,\,\sigma_y^{(N-1)}]\otimes\sigma_y^{(N)}+\\
&{}&I^{(1)}\otimes\bigl[A_1,\,\sum_{k=2}^{N-1}\,\sigma_z^{(k)}+\sum_{k=2}^{N-2}(\sigma_x^{(k)}\otimes\sigma_x^{k+1)}+\sigma_y^{(k)}\otimes\sigma_y^{(k+1)})\bigr]\otimes I^{(N)}=0.
\end{eqnarray*}
Hence we obtain the following
$$[A_1,\,\sigma_x^{(2)}]=[A_1,\,\sigma_y^{(2)}]=[A_1,\,\sigma_x^{(N-1)}]=[A_1,\,\sigma_y^{(N-1)}]=0.$$
This implies that
$$A_1=I^{(2)}\otimes A_2\otimes I^{(N-1)},$$
where $A_2$ is an operator on $\otimes_{k=3}^{N-2}\,\C^2$.

Repeating this argument until one arrives that $A=\lambda I$. 
Thus we obtain
$$\{H_S,\,\sigma_+^{(1)},\;\sigma_-^{(1)},\;\sigma_+^{(N)},\;\sigma_-^{(N)}\}'=\C I.$$
Therefore by Theorem \ref{a} we can conclude.
\end{proof}
\subsection{Approach to equilibrium}

The aim of this section is to prove that the quantum dynamical semigroup associated to the spin chain coupled to two heat baths has the property of approach to equilibrium.

The following theorem is introduced in [B]. 
\begin{theo}\label{th}
Let $L$ be a generator of a norm continuous quantum dynamical semigroup $(\Theta_t)_t$ on $\Bb(\Kk)$ which has the form given in (\ref{A}) and where the number of induces $j$ is finite. Assume that the following hypothesis hold:
\begin{enumerate}
\item[i)]  There exists a stationary state $\rho$ for the quantum dynamical semigroup $(\Theta_t^*)_t$,
\item[ii)] The linear span of all $V_j$ is self-adjoint,
\item[iii)] If $A\in\Bb(\Kk)$ such that $\Theta_t(A^*A)=(\Theta_tA^*)(\Theta_tA),$ for all $t\geq0$, then $A=\C I.$
\end{enumerate}
Then the state $\rho$ is faithful and the quantum dynamical semigroup $(\Theta_t^*)_t$ has the property of approach to equilibrium, that is
$$\lim_{t\rightarrow\infty}\Tr(\Theta_t^*\xi A)=\Tr(\rho A),\;\mbox{ for all normal state }\,\xi\,\mbox{ and }\;\mbox{ for all }\;A\in\Bb(\Kk).$$
\end{theo}

Under the hypothesis of the above theorem, $\rho$ is the unique stationary state for the quantum dynamical semigroup $(\Theta_t^*)_t$. In fact, let us consider an element $A$ in $\Bb(\Kk)$ such that $[H,A]=[V_j,A]=[V_j^*,A]=0$ for all $j$. Thus from hypothesis ii), $[V_j,A]=0$ implies that $[V_j^*,A]=0$. Hence we obtain
$$L(A)=L(A^*)=L(A^*A)=0.$$ 
It follows that $\Theta_tA^*=A^*$, $\Theta_tA=A$ and $\Theta_t(A^*A)=A^*A$ for all $t\geq0$. Then we get $\Theta_t(A^*A)=(\Theta_tA^*)(\Theta_tA)$, for all $t\geq0$. Finally, from the hypothesis iii), we have $A=\lambda I$. Note that $\rho$ is a faithful state. Therefore by Theorem \ref{a} we can conclude. 

Now as a corollary of Theorem \ref{th}, we prove the following.
\begin{theo}\label{equ}
The quantum dynamical semigroup $\{T_t^*=e^{t\mathcal{L}^*},\;t\in\R_+\}$ associated to the spin chain coupled to two heat baths at inverse temperatures $\beta$ and $\beta'$ has the property of approach to equilibrium to a unique stationary faithful state $\rho^{\beta,\,\beta'}$.
\end{theo} 
\begin{proof}
Note that dim $\Hh_S<\infty$. Then there exists a stationary state for the quantum dynamical semigroup $(T_t^*)_t$. This implies that assumption i) of the above theorem is satisfied. Moreover, it is clear that the linear span $\{\sigma_-^{(1)},\sigma_+^{(1)},\sigma_-^{(N)},\sigma_+^{(N)}\}$ is self-adjoint.

Now set $A\in\Bb(\Hh_S)$ such that 
\begin{equation}\label{iii)}
T_t(A^*A)=(T_tA^*)(T_tA),\;\,\forall \,t\geq0.
\end{equation}
By using the properties of the semigroups we can deduce that
\begin{eqnarray*}
T_s((T_tA)^*(T_tA))&=&T_s((T_tA^*)(T_tA))\\
                           &=&T_s(T_t(A^*A))\\
                           &=&T_{s+t}(A^*A)\\
                           &=&(T_{s+t}A^*)(T_{s+t}A)\\
                           &=&(T_s(T_tA)^*)(T_s(T_tA)),
 \end{eqnarray*}
for all $s\geq0$. Therefore for all $t\geq0$ the operator $T_tA$ satisfies relation (\ref{iii)}).

Note that by taking the derivative in (\ref{iii)}) with respect to $t$, we have
$$\mathcal{L}T_t(A^*A)=(\mathcal{L}T_tA^*)(T_tA)+(T_tA^*)(\mathcal{L}T_tA).$$
But we have $T_tA^*=(T_tA)^*$. Hence we obtain
\begin{equation}\label{B}
\mathcal{L}((T_tA^*)(T_tA))=(\mathcal{L}(T_tA)^*)(T_tA)+(T_tA^*)(\mathcal{L}T_tA),\,\;\forall t\,\geq0.
\end{equation}
In particular, for $t=0$ we have
\begin{equation}\label{C}
\mathcal{L}(A^*A)=(\mathcal{L}A^*)A+A^*(\mathcal{L}A).
\end{equation}
In the other hand we have the following 
\begin{eqnarray*}
\mathcal{L}(A^*A)-(\mathcal{L}A^*)A-A^*(\mathcal{L}A)&=&4\beta_0[\sigma_+^{(1)},\,A]^*\,[\sigma_+^{(1)},\,A]+4\beta_1[\sigma_-^{(1)},\,A]^*\,[\sigma_-^{(1)},\,A]\\
&+&4\beta_0'[\sigma_+^{(N)},\,A]^*\,[\sigma_+^{(N)},\,A]+4\beta_1'[\sigma_-^{(N)},\,A]^*\,[\sigma_-^{(N)},\,A].
\end{eqnarray*}
Hence if $A$ satisfies relation (\ref{iii)}), then the operator $A$ satisfies
$$A\in\{\sigma_-^{(1)},\sigma_+^{(1)},\sigma_-^{(N)},\sigma_+^{(N)}\}'.$$
 This gives $$A=I^{(1)}\otimes\widetilde{A}\otimes I^{(N)},$$
 where $\widetilde{A}$ is an operator on $\otimes_{k=2}^{N-1}\C^2$. Besides from relation (\ref{B}), the operator $T_t^*A$ satisfies also (\ref{C}). Therefore $T_tA$ has to be of the same form as $A$, that is 
$$T_tA=I^{(1)}\otimes\widetilde{S}_t\otimes I^{(N)},$$
 with $\widetilde{S}_t$ is an operator on $\otimes_{k=2}^{N-1}\C^2$. Furthermore by taking the derivative of $T_tA$ with respect to t at $t=0$ we obtain
\begin{eqnarray}\label{D}
\mathcal{L}(A)&=&i\,\bigl[H_S,\,A\bigr]\nonumber\\
              &=&i\,\sigma_x^{(1)}\otimes\,\bigl[\sigma_x^{(2)},\,\widetilde{A}\bigr]\otimes I^{(N)}+i\,\sigma_y^{(1)}\otimes\,\bigl[\sigma_y^{(2)},\,\widetilde{A}\bigr]\otimes I^{(N)}\nonumber\\
              &+&i \,I^{(1)}\otimes\,\Bigl[\sum_{k=2}^{N-2}(\sigma_x^{(k)}\otimes \sigma_x^{(k+1)}+\sigma_y^{(k)}\otimes \sigma_y^{(k+1)}),\,\widetilde{A}\Bigr]\otimes I^{(N)}\\
              &+&i \,I^{(1)}\otimes\,\bigl[\sigma_x^{(N-1)},\,\widetilde{A}\bigr]\otimes \sigma_x^{(N)}+i\,I^{(1)}\otimes\,\bigl[\sigma_y^{(N-1)},\,\widetilde{A}\bigr]\otimes \sigma_y^{(N)}\nonumber\\
              &=&I^{(1)}\otimes\widetilde{B}\otimes I^{(N)},\nonumber
\end{eqnarray}
where $\widetilde{B}=\frac{d}{dt}\widetilde{S}_t\bigl|_{t=0}$. Then from equality (\ref{D}) we have
$$[\sigma_x^{(2)},\,\widetilde{A}]=[\sigma_y^{(2)},\,\widetilde{A}]=[\sigma_x^{(N-1)},\,\widetilde{A}]=[\sigma_y^{(N-1)},\,\widetilde{A}]=0.$$
This implies that 
$$\widetilde{A}=I^{(2)}\otimes\widetilde{A}_1\otimes I^{(N-1)}$$
and
$$A=I^{(1)}\otimes I^{(2)}\otimes \widetilde{A}_1\otimes I^{(N-1)}\otimes I^{(N)}.$$
Note that $T_tA$ satisfies relation (\ref{iii)}) for all $t\geq0$. Hence by the same argument as before, $T_tA$ is written as
$$T_tA=I^{(1)}\otimes I^{(2)}\otimes \widetilde{R}_t\otimes I^{(N-1)}\otimes I^{(N)},$$
where $\widetilde{R}_t$ is an operator on $\otimes_{k=3}^{N-2}\C^2$. Repeating this reasoning until one obtains the result that is only possible if $A$ is a multiple of the identity. This ends the proof.
\end{proof}
\subsection{Local equilibrium states}
 Here we suppose that the spin chain is coupled to two heat baths at inverse temperatures $\beta$ and $\beta'$. Let us recall that there exists a unique stationary state $\rho^{\beta,\beta'}$ of the associated quantum dynamical semigroup. For $i\in\{1,...,N\}$, we denote by $\rho^{(i)}$ the local state associated to the $i$-th spin which is given by
$$\Tr(\rho^{(i)}A^{(i)})=\Tr(\rho^{\beta,\beta'}(I\otimes A^{(i)}\otimes I)),$$
where $A^{(i)}$ is an operator acting on the $i$-th copy of $\C^2$ in the chain $\otimes_{k=1}^{N}\C^2$.

 In this section, we treat the cases of the spin chain when it is made up of 2, 3 and 4 spins.\smallskip\\
- For $N=2$, we have
\begin{eqnarray*}
\rho^{\beta,\beta'}&=&(\frac{\rho_\beta+\rho_{\beta'}}{2})\otimes(\frac{\rho_\beta+\rho_{\beta'}}{2})-\frac{1}{8}(\beta_0-\beta_0')^2\sigma_z\otimes\sigma_z\\
& &+\frac{(\beta_0-\beta_0')}{4}\big[n_+\otimes n_--n_-\otimes n_+\big]+i\frac{(\beta_0-\beta_0'}{4}\big[\sigma_+\otimes\sigma_--\sigma_-\otimes\sigma_+\big].
\end{eqnarray*}
Therefore we get
\begin{eqnarray*}
\rho^{(1)}&=&\frac{\rho_\beta+\rho_{\beta'}}{2}+\frac{1}{2}\,(\frac{\rho_\beta-\rho_{\beta'}}{2}),\\
\rho^{(2)}&=&\frac{\rho_\beta+\rho_{\beta'}}{2}+\frac{1}{2}\,(\frac{\rho_{\beta'}-\rho_\beta}{2}).
\end{eqnarray*}

- For $N=3$, the equilibrium state $\rho^{\beta,\,\beta'}$ is given by
\begin{eqnarray*}
\rho^{\beta,\,\beta'}&=&(\frac{\rho_\beta+\rho_{\beta'}}{2})\otimes(\frac{\rho_\beta+\rho_{\beta'}}{2})\otimes(\frac{\rho_\beta+\rho_{\beta'}}{2})-\frac{3}{4}\,(\frac{\rho_\beta-\rho_{\beta'}}{2})\otimes(\frac{\rho_\beta+\rho_{\beta'}}{2})\otimes(\frac{\rho_\beta-\rho_{\beta'}}{2})\\
                   &+&\frac{3}{4}\,\big[(\frac{\rho_\beta-\rho_{\beta'}}{2})\otimes(\frac{\rho_\beta+\rho_{\beta'}}{2})\otimes(\frac{\rho_\beta+\rho_{\beta'}}{2})-(\frac{\rho_\beta+\rho_{\beta'}}{2})\otimes(\frac{\rho_\beta+\rho_{\beta'}}{2})\otimes(\frac{\rho_\beta-\rho_{\beta'}}{2})\big]\\
                  &+&\frac{\beta_0-\beta'_0}{8}\,\big[(\rho_\beta\otimes n_-\otimes n_+-\rho_\beta\otimes n_+\otimes n_-)+(n_-\otimes n_+\otimes\rho_{\beta'}-n_+\otimes n_-\otimes\rho_{\beta'})\big]\\
                  &+&i\,\frac{\beta_0-\beta'_0}{8}\,\big[\sigma_+\otimes\sigma_-\otimes(\frac{\rho_\beta+\rho_{\beta'}}{2})-\sigma_-\otimes\sigma_+\otimes(\frac{\rho_\beta+\rho_{\beta'}}{2})\big]\\
                  &+&i\,\frac{\beta_0-\beta'_0}{8}\,\big[(\frac{\rho_\beta+\rho_{\beta'}}{2})\otimes\sigma_+\otimes\sigma_--(\frac{\rho_\beta+\rho_{\beta'}}{2})\otimes\sigma_-\otimes\sigma_+\big]\\
                  &+&i\,\frac{\beta_0-\beta'_0}{8}\,\big[\rho_\beta\otimes\sigma_+\otimes\sigma_--\rho_\beta\otimes\sigma_-\otimes\sigma_+\big]\\
                  &+&i\,\frac{\beta_0-\beta'_0}{8}\,\big[\sigma_+\otimes\sigma_-\otimes\rho_{\beta'}-\sigma_-\otimes\sigma_+\otimes\rho_{\beta'}\big]\\
                 &-&\frac{(\beta_0-\beta_0')^2}{16}\,\big[\sigma_+\otimes I\otimes\sigma_-+\sigma_-\otimes I\otimes\sigma_+\big].
\end{eqnarray*}
Thus the local states are given by
\begin{eqnarray*}
\rho^{(1)}&=&\frac{\rho_\beta+\rho_{\beta'}}{2}+\frac{1}{2}\,(\frac{\rho_\beta-\rho_{\beta'}}{2}),\\
\rho^{(2)}&=&\frac{\rho_\beta+\rho_{\beta'}}{2}=\frac{\rho^{(1)}+\rho^{(3)}}{2},\\
\rho^{(3)}&=&\frac{\rho_\beta+\rho_{\beta'}}{2}+\frac{1}{2}\,(\frac{\rho_{\beta'}-\rho_\beta}{2}).
\end{eqnarray*}

- For $N=4$, we have
\begin{eqnarray*}
\rho^{\beta,\,\beta'}\!\!&=&\!\!\!(\frac{\rho_\beta+\rho_{\beta'}}{2})\otimes(\frac{\rho_\beta+\rho_{\beta'}}{2})\otimes(\frac{\rho_\beta+\rho_{\beta'}}{2})\otimes(\frac{\rho_\beta+\rho_{\beta'}}{2})\\
& &\!\!\!\!-\frac{7}{8}\,(\frac{\rho_\beta-\rho_{\beta'}}{2})\otimes(\frac{\rho_\beta+\rho_{\beta'}}{2})\otimes(\frac{\rho_\beta+\rho_{\beta'}}{2})\otimes(\frac{\rho_\beta-\rho_{\beta'}}{2})\\
& &+\frac{1}{2}\,\big[(\frac{\rho_\beta-\rho_{\beta'}}{2})\otimes(\frac{\rho_\beta+\rho_{\beta'}}{2})\otimes(\frac{\rho_\beta+\rho_{\beta'}}{2})\otimes(\frac{\rho_\beta+\rho_{\beta'}}{2})\\
& &-(\frac{\rho_\beta+\rho_{\beta'}}{2})\otimes(\frac{\rho_\beta+\rho_{\beta'}}{2})\otimes(\frac{\rho_\beta+\rho_{\beta'}}{2})\otimes(\frac{\rho_\beta-\rho_{\beta'}}{2})\big]\\
& &-\frac{1}{8}\,\rho_\beta\otimes(\frac{\rho_\beta-\rho_{\beta'}}{2})\otimes(\frac{\rho_\beta-\rho_{\beta'}}{2})\otimes\rho_{\beta'}\\
& &-\frac{(\beta_0-\beta_0')^2}{32}\,\big[n_-\otimes(\frac{\rho_\beta+\rho_{\beta'}}{2})\otimes I\otimes n_++n_+\otimes I\otimes(\frac{\rho_\beta+\rho_{\beta'}}{2})\otimes n_-\big]\\
& &-\frac{(\beta_0-\beta_0')^2}{16}\,\big[(\frac{\rho_\beta+\rho_{\beta'}}{2})\otimes n_+\otimes n_+\otimes n_-+(\frac{\rho_\beta+\rho_{\beta'}}{2})\otimes n_-\otimes n_-\otimes n_+\big]\\
& &-\frac{(\beta_0-\beta_0')^2}{16}\,\big[n_-\otimes n_+\otimes n_+\otimes(\frac{\rho_\beta+\rho_{\beta'}}{2})+n_+\otimes n_-\otimes n_-\otimes(\frac{\rho_\beta+\rho_{\beta'}}{2})\big]\\
& &+\frac{(\beta_0-\beta_0')^2}{16}\,\big[n_-\otimes(\frac{\rho_\beta+\rho_{\beta'}}{2})\otimes n_+\otimes n_-+n_+\otimes n_-\otimes(\frac{\rho_\beta+\rho_{\beta'}}{2})\otimes n_+\big]\\
& &+\frac{(\beta_0-\beta_0')^2}{32}\,\big[\rho_\beta\otimes n_-\otimes n_+\otimes n_-+n_+\otimes n_+\otimes n_-\otimes\rho_\beta\big]\\
& &+\frac{(\beta_0-\beta_0')^2}{32}\,\big[\rho_{\beta'}\otimes n_-\otimes n_+\otimes n_++n_-\otimes n_+\otimes n_-\otimes\rho_{\beta'}\big]\\
& &+\frac{(\beta_0-\beta_0')^2}{32}\,(\beta_0+\beta_0')\,n_+\otimes n_+\otimes n_-\otimes n_+\\
& &+\frac{(\beta_0-\beta_0')^2}{32}\,(\beta_1+\beta_1')\,n_-\otimes n_+\otimes n_-\otimes n_-\\
& &+\frac{3(\beta_0-\beta_0')^2}{32}\,\big[n_+\otimes n_-\otimes n_-\otimes n_++n_-\otimes n_+\otimes n_+\otimes n_-\big]\\
& &-\frac{(\beta_0-\beta_0')^2}{16}\,\big[n_+\otimes n_+\otimes n_-\otimes n_-+n_-\otimes n_-\otimes n_+\otimes n_+\big]\\  
& &+i\,\frac{(\beta_0-\beta_0')}{32}\,(\frac{\rho_\beta+\rho_{\beta'}}{2})\otimes(\frac{\rho_\beta+\rho_{\beta'}}{2})\otimes\sigma_+\otimes\sigma_-\\
& &-i\,\frac{(\beta_0-\beta_0')}{32}\,(\frac{\rho_\beta+\rho_{\beta'}}{2})\otimes(\frac{\rho_\beta+\rho_{\beta'}}{2})\otimes\sigma_-\otimes\sigma_+\\
& &+i\,\frac{(\beta_0-\beta_0')}{32}\,\sigma_+\otimes\sigma_-\otimes(\frac{\rho_\beta+\rho_{\beta'}}{2})\otimes(\frac{\rho_\beta+\rho_{\beta'}}{2})\\
& &-i\,\frac{(\beta_0-\beta_0')}{32}\,\sigma_-\otimes\sigma_+\otimes(\frac{\rho_\beta+\rho_{\beta'}}{2})\otimes(\frac{\rho_\beta+\rho_{\beta'}}{2})\\
& &+i\,\frac{(\beta_0-\beta_0')}{16}\,\big[\rho_\beta\otimes\rho_{\beta'}\otimes\sigma_+\otimes\sigma_--\rho_\beta\otimes\rho_{\beta'}\otimes\sigma_-\otimes\sigma_+\big]\\
& &+i\,\frac{(\beta_0-\beta_0')}{16}\,\big[\sigma_+\otimes\sigma_-\otimes\rho_\beta\otimes\rho_{\beta'}-\sigma_-\otimes\sigma_+\otimes\rho_\beta\otimes\rho_{\beta'}\big]\\
& &+i\,\frac{(\beta_0-\beta_0')}{8}\,\big[\sigma_+\otimes\sigma_-\otimes(\frac{\rho_\beta+\rho_{\beta'}}{2})\otimes\rho_{\beta'}-\sigma_-\otimes\sigma_+\otimes(\frac{\rho_\beta+\rho_{\beta'}}{2})\otimes\rho_{\beta'}\big]\\
& &+i\,\frac{(\beta_0-\beta_0')}{8}\,\big[\rho_\beta\otimes(\frac{\rho_\beta+\rho_{\beta'}}{2})\otimes\sigma_+\otimes\sigma_--\rho_\beta\otimes(\frac{\rho_\beta+\rho_{\beta'}}{2})\otimes\sigma_-\otimes\sigma_+\big]\\
\end{eqnarray*}
\begin{eqnarray*}
& &+i\,\frac{3(\beta_0-\beta_0')}{32}\,(\frac{\rho_\beta+\rho_{\beta'}}{2})\otimes\sigma_+\otimes\sigma_-\otimes(\frac{\rho_\beta+\rho_{\beta'}}{2})\\
& &-i\,\frac{3(\beta_0-\beta_0')}{32}\,(\frac{\rho_\beta+\rho_{\beta'}}{2})\otimes\sigma_-\otimes\sigma_+\otimes(\frac{\rho_\beta+\rho_{\beta'}}{2})\big]\\
& &-i\,\frac{(\beta_0-\beta_0')^2}{32}\,\big[n_+\otimes n_-\sigma_+\otimes\sigma_--n_+\otimes n_-\otimes\sigma_-\otimes\sigma_+\big]\\
& &+i\,\frac{(\beta_0-\beta_0')^2}{32}\,\big[n_-\otimes n_+\otimes\sigma_+\otimes\sigma_--n_-\otimes n_+\otimes\sigma_-\otimes\sigma_+\big]\\
& &-i\,\frac{(\beta_0-\beta_0')^2}{32}\,\big[\sigma_+\otimes\sigma_-\otimes n_+\otimes n_--\sigma_-\otimes\sigma_+\otimes n_+\otimes n_-\big]\\
& &-i\,\frac{(\beta_0-\beta_0')^2}{32}\,\big[\sigma_+\otimes\sigma_-\otimes n_-\otimes n_+-\sigma_-\otimes\sigma_+\otimes n_-\otimes n_+\big]\\
& &+i\,\frac{(\beta_0-\beta_0')^2}{32}\,\big[n_+\otimes\sigma_+\otimes\sigma_-\otimes n_--n_+\otimes\sigma_-\otimes\sigma_+\otimes n_-\big]\\
& &+i\,\frac{(\beta_0-\beta_0')^2}{32}\,\big[n_-\otimes\sigma_+\otimes\sigma_-\otimes n_+-n_-\otimes\sigma_-\otimes\sigma_+\otimes n_+\big]\\
& &+\big(\frac{1}{64}\,(\beta_0+\beta_0')^2\,-\frac{1}{16}\,\beta_0^2\,\big)\,(\beta_0-\beta_0')\,\big[\sigma_z\otimes\sigma_+\otimes I\otimes\sigma_-+\sigma_z\otimes\sigma_-\otimes I\otimes\sigma_+\big]\\
& &-\big(\frac{1}{64}\,(\beta_0+\beta_0')^2-\frac{1}{16}\,\beta_0'^2\big)\,(\beta_0-\beta_0')\,\big[\sigma_+\otimes I\otimes\sigma_-\otimes\sigma_z+\sigma_-\otimes I\otimes\sigma_+\otimes\sigma_z\big]\\
& &-\frac{(\beta_0-\beta_0')^2}{16}\,\big[n_-\otimes\sigma_+\otimes I\otimes\sigma_-+n_-\otimes\sigma_-\otimes I\otimes\sigma_+\big]\\
& &-\frac{(\beta_0-\beta_0')^2}{16}\,\big[\sigma_+\otimes I\otimes\sigma_-\otimes n_-+\sigma_-\otimes I\otimes\sigma_+\otimes n_-\big]\\
& &+\frac{(\beta_0-\beta_0')^3}{64}\,\big[I\otimes\sigma_+\otimes\sigma_-\otimes I+I\otimes\sigma_-\otimes\sigma_+\otimes I\big]\\
& &+\frac{(\beta_0-\beta_0')^2}{16}\,\big[\sigma_+\otimes\sigma_-\otimes\sigma_+\otimes\sigma_-+\sigma_-\otimes\sigma_+\otimes\sigma_-\otimes\sigma_+\big]\\
& &-\frac{(\beta_0-\beta_0')^2}{16}\,\big[\sigma_+\otimes\sigma_-\otimes\sigma_-\otimes\sigma_++\sigma_-\otimes\sigma_+\otimes\sigma_+\otimes\sigma_-\big].
\end{eqnarray*}
Hence we obtain
\begin{eqnarray*}
\rho^{(1)}&=&\frac{\rho_\beta+\rho_{\beta'}}{2}+\frac{1}{2}\,(\frac{\rho_\beta-\rho_{\beta'}}{2}),\\
\rho^{(2)}&=&\rho^{(3)}=\frac{\rho_\beta+\rho_{\beta'}}{2}=\frac{\rho^{(1)}+\rho^{(4)}}{2},\\
\rho^{(4)}&=&\frac{\rho_\beta+\rho_{\beta'}}{2}+\frac{1}{2}\,(\frac{\rho_{\beta'}-\rho_\beta}{2}).
\end{eqnarray*}  

Note that for $N\geq5$, it is very hard to find the stationary state $\rho^{\beta,\,\beta'}$. Moreover, from the computation done in the cases
of the spin chain when it is made up of 2, 3 and 4 atoms, we see that when $N$ increase, the number of the off-diagonal terms increase quickly in the explicit form of the matrix of the state $\rho^{\beta,\,\beta'}$ in the canonical basis of $\Hh_S$. Besides the off-diagonal terms do not contribute to the calculation of the partial trace at any site. However the form of the diagonal terms given in the cases $N=$2, 3 and 4 are similar enough that we believe that the following conjecture is true: For $N\geq5$, the local states are given by 
\begin{eqnarray*}
\rho^{(1)}&=&\frac{\rho_\beta+\rho_{\beta'}}{2}+\frac{1}{2}\,(\frac{\rho_\beta-\rho_{\beta'}}{2}),\\
\rho^{(2)}&=&...=\rho^{(N-1)}=\frac{\rho_\beta+\rho_{\beta'}}{2}=\frac{\rho^{(1)}+\rho^{(N)}}{2},\\
\rho^{(N)}&=&\frac{\rho_\beta+\rho_{\beta'}}{2}+\frac{1}{2}\,(\frac{\rho_{\beta'}-\rho_\beta}{2}).
\end{eqnarray*}
\subsection{Entropy production}
In this section, we treat the case of the spin chain coupled to two heat baths at the same temperature $\beta^{-1}$. Let us recall that from Theorem \ref{same-temperature}, there exists a unique stationary faithful state $\rho^\beta$ for the associated quantum master equation in the case of same temperature $\beta=\beta'$. The definition of entropy production, that we give here, is taken from [SL].

Let $\rho$ be a state on $\Hh_S$ and set $\rho(t)=e^{t\mathcal{L}^*}(\rho)$. Then the relative entropy of $\rho$ with respect to $\rho^\beta$ is defined by
$$S(\rho(t)|\rho^\beta)=\Tr(\rho(t)(\log\rho^\beta-\log\rho(t))).$$
Hence the entropy production is given by
\begin{eqnarray*}
\sigma(\rho)&=&-\frac{d}{dt}\,S(\rho(t)|\rho^\beta)\bigl|_{t=0}\\
            &=&\Tr(\mathcal{L}^*(\rho)(\log\rho^\beta-\log\rho)),
\end{eqnarray*}
 where $\Tr(\mathcal{L}^*(\rho)\log\rho)$ is given as
\begin{eqnarray*}
\Tr(\mathcal{L}^*(\rho)\log\rho)&=&\sum_j\,\langle\Psi_j,\,\mathcal{L}^*(\rho)\Psi_j\rangle\,\log\rho_j,\\
\langle\Psi_j,\,\mathcal{L}^*(\rho)\Psi_j\rangle\log\rho_j&=&
\left\{
\begin{array}{lcc}
-\infty  &\mbox{ if }&\langle\Psi_j,\,\mathcal{L}^*(\rho)\Psi_j\rangle\neq0 \mbox{ and }\,\rho_j=0\\
0  &\mbox{ if }& \langle\Psi_j,\,\mathcal{L}^*(\rho)\Psi_j\rangle=0.
\end{array}
\right.
\end{eqnarray*}
\begin{theo}
The entropy production associated to the spin chain coupled to two heat baths at the same inverse temperature $\beta$ is written as
$$\sigma(\rho)=4\beta_0\big[\sum_{j,\,k}\big[|\langle\Psi_k,\,\sigma_+^{(1)}\Psi_j\rangle|^2+|\langle\Psi_k,\,\sigma_+^{(N)}\Psi_j\rangle|^2]\,(e^{2\beta}\rho_k-\rho_j)(\log\rho_k-\log\rho_j+2\beta)\big],$$
where $\rho=\sum_j\,\rho_j|\Psi_j\rangle\langle\Psi_j|$ is the spectral decomposition of $\rho$.
\end{theo}
\begin{proof}
Note that 
$$\mathcal{L}^*=\mathcal{L}^*_h+\mathcal{L}^*_d,$$
where $\mathcal{L}^*_h$ is the Hamiltonian part of $\mathcal{L}^*$ and $\mathcal{L}^*_d=\mathcal{L}_d^{*(1)}+\mathcal{L}_d^{*(N)}$ is its dissipative part with
\begin{eqnarray*}
\mathcal{L}_d^{*(1)}(\rho)&=&2\beta_0\,[2\sigma_+^{(1)}\rho\,\sigma_-^{(1)}-\{n_-^{(1)},\,\rho\}]\\
                   &+&2\beta_1\,[2\sigma_-^{(1)}\rho\,\sigma_+^{(1)}-\{n_+^{(1)},\,\rho\}],\\
\mathcal{L}_d^{*(N)}(\rho)&=&2\beta_0'\,[2\sigma_+^{(N)}\rho\,\sigma_-^{(N)}-\{n_-^{(N)},\,\rho\}]\\
                   &+&2\beta_1'\,[2\sigma_-^{(N)}\rho\,\sigma_+^{(N)}-\{n_+^{(N)},\,\rho\}].
\end{eqnarray*} 

Put
$$H^{(S)}=\sum_{k=1}^N\,\sigma_z^{(k)}.$$
It is easy to show that the equilibrium state $\rho^\beta$ is given by
$$\rho^\beta=\frac{1}{Z}\,e^{-\beta H^{(S)}},$$
where $Z=\Tr(e^{-\beta H^{(S)}}).$ Thus we obtain $\log\rho^\beta=-\beta H^{(S)}-\log Z$. On the other hand, a straightforward computation shows that
\begin{eqnarray*}
\Tr([H_S,\,\rho]\log\rho)=\Tr(H_S[\rho,\,\log\rho])=0
\end{eqnarray*}
and
\begin{eqnarray*}
\Tr([H_S,\,\rho]\log\rho^\beta)=-\beta\Tr([H^{(S)},\,H_S]\rho)=0.
\end{eqnarray*}
Therefore we get 
$$\Tr(\mathcal{L}^*_h(\rho(\log\rho^\beta-\log\rho)))=0.$$ 
This gives
\begin{equation}\label{E}
\sigma(\rho)=\sigma_1(\rho)+\sigma_N(\rho)=-\Tr(\mathcal{L}^*_d(\rho)\log\rho)-\beta\,\Tr(\mathcal{L}^*_d(\rho)H^{(S)}),
\end{equation}
where $\sigma_i(\rho)=-\Tr(\mathcal{L}^{*(i)}_d(\rho)\log\rho)-\beta\,\Tr(\mathcal{L}^{*(i)}_d(\rho)H^{(S)})$ with $i=1,\;N.$

Now let us compute the terms of the second member in (\ref{E}). We have
\begin{eqnarray*}
\Tr(\mathcal{L}^{*(1)}_d(\rho)\log\rho)&=&4\beta_0\Big[\sum_{j,\,k}\,\langle\Psi_j,\,\sigma_-^{(1)}\,\Psi_k\rangle\,\langle\Psi_k,\,\sigma_+^{(1)}\,\Psi_j\rangle\,\rho_j\log\rho_k\\
         &-&\sum_j\,\langle\Psi_j,\,n_-^{(1)}\Psi_j\rangle\,\rho_j\log\rho_j\Big]\\
         &+&4\beta_1\,\Big[\sum_{j,\,k}\,\langle\Psi_j,\,\sigma_+^{(1)}\Psi_k\rangle\,\langle\Psi_k,\,\sigma_-^{(1)}\Psi_j\rangle\,\rho_j\log\rho_k\\
         &-&\sum_j\,\langle\Psi_j,\,n_+^{(1)}\Psi_j\rangle\,\rho_j\log\rho_j\Big],\\
\Tr(\mathcal{L}^{*(1)}_d(\rho)H^{(S)})&=&8\beta_0\,\sum_j\,\langle\Psi_j,\,n_-^{(1)}\Psi_j\rangle\,\rho_j-8\beta_1\,\sum_j\,\langle\Psi_j,\,n_+^{(1)}\,\Psi_j\rangle\,\rho_j.
\end{eqnarray*}
Note that
\begin{eqnarray*}
\langle\Psi_j,\,n_+^{(1)}\Psi_j\rangle&=&||\sigma_+^{(1)}\Psi_j||^2=\sum_k\,|\langle\Psi_k,\,\sigma_+^{(1)}\Psi_j\rangle|^2,\\
\langle\Psi_j,\,n_-^{(1)}\Psi_j\rangle&=&||\sigma_-^{(1)}\Psi_j||^2=\sum_k\,|\langle\Psi_j,\,\sigma_+^{(1)}\Psi_k\rangle|^2.
\end{eqnarray*}
Therefore we obtain
\begin{eqnarray}\label{F}
 \sigma_1(\rho)&=&4\beta_0\Bigl[\sum_{j,\,k}|\langle\Psi_k,\,\sigma_+^{(1)}\Psi_j\rangle|^2\,\rho_j(\log\rho_j-\log\rho_k-2\beta)\Bigr]\nonumber\\
               &+&4\beta_1\Big[\sum_{j,\,k}\,|\langle\Psi_j,\,\sigma_+^{(1)}\Psi_k\rangle|^2\rho_j(\log\rho_j-\log\rho_k+2\beta)\Big].
\end{eqnarray}
If we substitute $\beta_1$ by $e^{2\beta}\beta_0$ in (\ref{F}), then we get 
$$\sigma_1(\rho)=4\beta_0\Bigl[\sum_{j,\,k}|\langle\Psi_k,\,\sigma_+^{(1)}\Psi_j\rangle|^2\,(e^{2\beta}\rho_k-\rho_j)(\log\rho_k-\log\rho_j+2\beta)\Bigr].$$
In the same way, we prove that
$$\sigma_N(\rho)=4\beta_0\Bigl[\sum_{j,\,k}|\langle\Psi_k,\,\sigma_+^{(N)}\Psi_j\rangle|^2\,(e^{2\beta}\rho_k-\rho_j)(\log\rho_k-\log\rho_j+2\beta)\Bigr].$$
This ends the proof of the above theorem.
\end{proof}

\begin{rmq}:
Note that as a corollary of the above theorem, we have $\sigma(\rho)\geq0$ for any density matrix $\rho$.
\end{rmq} 
\subsection{Quantum detailed balance condition}
In this section we suppose that the spin chain is coupled to two heat baths at same inverse temperature $\beta$. Let us recall that 
$$\rho^\beta=\otimes_{k=1}^N\rho_\beta$$ 
is the only stationary faithful state of the quantum dynamical semigroup $(T_t^*)_t$.  

The following definition is introduced in [AL].
\begin{defi}
Let $\Theta$ be a generator of a quantum dynamical semigroup written as
$$ \Theta=-i\,[H,.]+\,\Theta_0,$$
where $H$ is a self-adjoint operator. We say that $\Theta$ satisfies a quantum detailed balance condition with respect to a stationary state $\rho$ if
\begin{enumerate}
\item[i)] $[H,\,\rho]=0,$
\item[ii)] $\langle\Theta_0(A),\,B\rangle_\rho=\langle A,\,\Theta_0(B)\rangle_\rho,$ for all $A,\,B\in D(\Theta_0),$
\item[{}] with $\langle A,\,B\rangle_\rho=\Tr(\rho A^*B).$
\end{enumerate}
\end{defi}
Now we prove the following.
\begin{theo}
The generator $\mathcal{L}^*$ of the quantum dynamical semigroup of the spin chain coupled to two heat baths at same inverse temperature $\beta$ satisfies a quantum detailed balance condition with respect to the stationary state $\rho^\beta$.
\end{theo}
\begin{proof}
Note that 
$$\mathcal{L}^*=-i\,[H_S,\,.]+\mathcal{L}_d^*,$$
where $\mathcal{L}_d^*$ is the dissipative part. On the other hand, we have $[H_S,\,\rho^\beta]=0$. This proves that assumption i) of the above definition is satisfied. Furthermore it is easy to show that $\mathcal{L}_d^*$ is a self-adjoint operator with respect to the scalar product $\langle A,\,B\rangle_{\rho^\beta}$. Thus the above theorem holds.
\end{proof}
\section{Spin chain coupled to several heat baths}
Let us consider now a spin chain ($N$ spins) coupled to $r$ heat baths at inverse temperatures $\beta^{(k_1)},\beta^{(k_2)},...,\beta^{(k_r)}$, where $2\leq r\leq N$ for all $j=1,...,r$ and $k_j$ is the $k_j$-th site of the chain $\otimes_{i=1}^N\C^2$. The quantum repeated interaction Hamiltonian is given by
$$H=H_S\otimes I+I\otimes H_R+\frac{1}{\sqrt{h}}\,\sum_{j=1}^r\,(\sigma_x^{(k_j)}\otimes\sigma_x^{(k_j)}+\sigma_y^{(k_j)}\otimes\sigma_y^{(k_j)}).$$
Therefore we prove in the same way as in subsection 2.2 that the associated Lindblad generator is written as
\begin{eqnarray*}
\mathcal{L}(X)=i\,[H_S,\,X]&+&2\beta_0^{(k_1)}\,[2\sigma_-^{(k_1)}X\sigma_+^{(k_1)}-\{n_-^{(k_1)},\,X\}]\\
                           &+&2\beta_1^{(k_1)}\,[2\sigma_+^{(k_1)}X\sigma_-^{(k_1)}-\{n_+^{(k_1)},\,X\}]\\
                           &+&2\beta_0^{(k_2)}\,[2\sigma_-^{(k_2)}X\sigma_+^{(k_2)}-\{n_-^{(k_2)},\,X\}]\\
                           &+&2\beta_1^{(k_2)}\,[2\sigma_+^{(k_2)}X\sigma_-^{(k_2)}-\{n_+^{(k_2)},\,X\}]\\
                           &.&                         \\
                           &.&                         \\                         
                           &+&2\beta_0^{(k_r)}\,[2\sigma_-^{(k_r)}X\sigma_+^{(k_r)}-\{n_-^{(k_r)},\,X\}]\\
                           &+&2\beta_1^{(k_r)}\,[2\sigma_+^{(k_r)}X\sigma_-^{(k_r)}-\{n_+^{(k_r)},\,X\}],
\end{eqnarray*}
for all $X\in\Bb(\Hh_S).$

Actually we prove the following.
\begin{theo}
If $\beta^{(k_1)}=\beta^{(k_2)}=...=\beta^{(k_r)}=\beta$, then $\rho^\beta=\otimes_{i=1}^N\rho_\beta$ is the unique stationary state of the quantum dynamical semigroup $T_t^*=e^{t\mathcal{L}^*}$.
\end{theo}
\begin{proof}
The proof of this theorem is similar as the one of Theorem \ref{same-temperature}.
\end{proof}

The following theorem can be proved in the same way as Theorem \ref{equ}.
\begin{theo}
The quantum dynamical semigroup $T_t^*=e^{t\mathcal{L}^*}$ associated to the spin chain coupled to $r$ heat baths at inverse temperatures $\beta^{(k_1)},\;\beta^{(k_2)},\;,...,\beta^{(k_r)}$ has the property of approach to equilibrium to a unique stationary state.
\end{theo}

Now our purpose is to study the associated entropy production. Here we suppose that $\beta^{(k_1)}=\beta^{(k_2)}=...=\beta^{(k_r)}=\beta$. Let us put
$$\sigma_{k_i}(\rho)=4\beta_0\Bigl[\sum_{j,\,m}\big[|\langle\Psi_m,\,\sigma_+^{(k_i)}\Psi_j\rangle|^2(e^{2\beta}\rho_m-\rho_j)(\log\rho_m-\log\rho_j+2\beta)\Bigr].$$
Then it is easy to show that $\sigma_i(\rho)$ is the entropy production of the spin chain coupled to the $i$-th heat bath at the $k_i$-th spin. 
\begin{theo}
If $\beta^{(k_1)},\;\beta^{(k_2)},\;,...,\beta^{(k_r)}$, then the entropy production of the spin chain coupled to $r$ heat baths is given by
\begin{equation}\label{cor}
\sigma(\rho)=\sum_{i=1}^r\,\sigma_{k_i}(\rho).
\end{equation} 
\end{theo}
\begin{proof}
We prove in the same way as in Theorem 3.5 that the entropy production of the spin chain coupled to $r$ heat baths is given by
\begin{eqnarray*}
\sigma(\rho)&=&-\Tr(\mathcal{L}^*_d(\rho)\log\rho)-\beta\,\Tr(\mathcal{L}_d^*(\rho)H^{(S)})\\
            &=&\sum_{i=1}^r[-\Tr(\mathcal{L}^{*(k_i)}_d(\rho)\log\rho)-\beta\,\Tr(\mathcal{L}_d^{*(k_i)}(\rho)H^{(S)})]\\
            &=&\sum_{i=1}^r\sigma_{k_i}(\rho),
\end{eqnarray*}
where 
\begin{eqnarray*}
\mathcal{L}^{*(k_i)}_d(\rho)&=&2\beta_0^{(k_i)}\,[2\sigma_-^{(k_i)}X\sigma_+^{(k_i)}-\{n_-^{(k_i)},\,X\}]\\
                           &+&2\beta_1^{(k_i)}\,[2\sigma_+^{(k_i)}X\sigma_-^{(k_i)}-\{n_+^{(k_i)},\,X\}].
\end{eqnarray*}
But from the proof of Theorem 3.5, we have
$$\sigma_{k_i}(\rho)=4\beta_0\Big[\sum_{j,\,m}\big[|\langle\Psi_m,\,\sigma_+^{(k_i)}\Psi_j\rangle|^2(e^{2\beta}\rho_m-\rho_j)(\log\rho_m-\log\rho_j+2\beta)\Big].$$
This proves our theorem.
\end{proof}

As a consequence of (\ref{cor}), we conclude the following.
\begin{corollary}
We have $\sigma(\rho)\geq 0$ for any density matrix $\rho$. Moreover, $\sigma(\rho)=0$ if and only if $\sigma_{k_i}(\rho)=0$ for all $i=1,...,r$.
\end{corollary}

Note that if $\beta^{(k_1)}=\beta^{(k_2)}=...=\beta^{(k_r)}=\beta$, then it is easy to show that the quantum dynamical semigroup of the spin chain coupled to $r$ heat baths satisfies a quantum detailed balance condition with respect to the stationary state $\rho^\beta$.

{\small

\end{document}